\documentclass[letter,12pt]{article}
\usepackage{amsfonts}
\usepackage{amssymb}
\usepackage{amsmath}
\usepackage{color,graphicx}
\usepackage{cite}
\usepackage{subfigure}
\usepackage{float}
\usepackage{pstricks-add}
\newtheorem{theo}{Theorem}
\newtheorem{prop}[theo]{Proposition}
\newtheorem{lema}[theo]{Lemma}

\newtheorem{obs}[theo]{Remark}
\newcommand{\R}{\mathbb{R}}

\newcommand{\C}{\mathbb{C}}
\newcommand{\N}{\mathbb{N}}
\newcommand{\eps}{\varepsilon}

\title{Old and new about equidistant sets and   generalized  conics}
\author{Mario Ponce \and Patricio Santib\'a\~nez}
\date{}
\begin{document}
\maketitle
\begin{abstract}
This article is devoted to the study of classical and new results concerning equidistant sets, both from the topological and metric point of view. We start with a review of the most interesting known facts about these sets in the euclidean space and then we prove that equidistant sets vary continuously with  their focal sets.    In the second part we propose a viewpoint in which equidistant sets can be thought of as natural generalization for conics. Along these lines, we show that many geometric features of classical conics can be retrieved in more general equidistant sets. In the Appendix  we prove a shadowing property of equidistant sets and provide sharp estimates. This result should be of interest for computer simulations.

\end{abstract}
\section{Introduction}
The set of points that are equidistant from two given  sets in the plane appears na\-tu\-ra\-lly  in many classical geometric situations. Namely, the classical conics, defined as the level set of a degree $2$ real polynomial equation, can always be realized  as the equidistant set to two circles (see section \ref{conics_as_midsets}).  The significance of conics for the development of Mathematics is incontestable.  Each new progress in their  study has represented a breakthrough: the determination of the bounded area by Archimides, their conception  as plane curves by Apollonius, their occurrence as solutions  for movement equations by Kepler,  the development of projective and analytic geometry by Desargues and Descartes, etc.\\

In other, let us say, less academic field we find equidistant  sets as conventionally defined frontiers in territorial domain controversies: for instance, the United Nations Convention on the Law of the Sea (Article 15) establishes that, in absence of any previous agreement,   the delimitation of the territorial sea between countries occurs exactly on the {\it median line every point of which is equidistant of the nearest points} to each country. The significance of this last situation  stresses the necessity of understanding the geometric structure of equidistant sets. The study of equidistant sets, other than conics, arose many decades ago principally with the works of  Wilker \cite{WILK75} and Loveland \cite{LOVE76} (in section \ref{survey} we review their main contributions concerning topological properties of equidistant sets).  \\

 In the literature we can find many generalizations of conics. For instance, Gro\ss { }and Strempel \cite{GRST98} start from the the usual definition of conics as the set of points in the plane that have  a constant weighted sum of distances to two points (the {\it focal points}) and gave a generalization by allowing more than two focal points, weights other than $\pm 1$ and point sets in higher dimensions. Recently, Vincze and Nagy \cite{VINA2010} proposed  that a generalized conic  is a set of points with the same average distance from a point set $\Gamma\subset \R^n$. In this work we present equidistant sets as a natural ge\-ne\-ra\-li\-za\-ti\-on of conics making use of the fact that classical conics are equidistant sets to plane circles ({\it focal sets}). Admitting more complicated focal sets we obtain more complicated equidistant sets ({\it generalized conics}).  Our purpose is to show that these generalizations share many geometrical features  with their classical ancestors. For instance,  we show that  the equidistant set to two disjoint connected compact sets  looks like the branch of a hyperbola in the sense that near the infinity it  is asymptotic to two rays.  We also discuss possible generalizations of ellipses and parabolas. Additionally we propose some further research directions.\\

Computational simulations constitute a useful tool for treating equidistant sets. This task faces two main theoretic  issues. The first one deals with the fact that a computer manipulate only discrete approximations of plane sets (limited by the memory or screen resolution). In this direction our result about the Hausdorff continuity of equidistant sets (cf. Theorem \ref{hausdorff}) is central since it implies that the equidistant set computed by the machine approaches (by increasing the screen resolution or dedicated memory) to the genuine equidistant set.  The second serious computational issue has to do with the following: remember that the plane is represented on the screen by a finite array of points (pixels). In order to determine precisely what points belong to the equidistant set, one needs to compute the respective distances to the underlying sets (focal sets) and then to decide whether this difference is equal to zero or not. Strictly speaking, a pixel belongs to the equidistant set if and only if this difference vanishes. Nevertheless, due to  to the discrete character of this  computer screen plane, this absolute zero is virtually impossible. Thus, in order to obtain a good picture of the equidistant set we need to introduce a more tolerant criteria. The general situation is as follows: one finds a pixel for which the difference of the distances to the focal sets  is very small  and we ask wether or not  this means that we can assume the presence of a point of the   equidistant set  inside the region represented by this pixel. Theorem \ref{shadow} (the {\it Shadowing property}, see the Appendix)  provides us with a useful criteria, since it gives a sharp bound on the distance from a quasi-equidistant point to an authentic equidistant point.

\paragraph{Definitions and notations. }
We consider $\R^n$ endowed with the classical euclidean distance $dist(\cdot, \cdot)$. One easily extends the definition to admit the distance between a point and a set. Given two non-empty sets $A, B\subset \R^n$ we define the {\it equidistant set} to $A$ and $B$ as
\[
\{A=B\}:=\{x\in \R^n \ : \ dist(x, A)=dist(x, B)\}.
\]
This notation is due to Wilker \cite{WILK75}. We also utilize the word {\it midset} as proposed  by Loveland \cite{LOVE76}. We say that $A$ and $B$ are the {\it focal sets} of the midset $\{A=B\}$.
\\

\noindent
For $x\in \R^n$ we write $\mathcal{P}_x(A)=\{p\in A \ : \ dist(x, A)=dist(x, p)\}$ the set of {\it foot points from $x$ to $A$}. \\

\noindent
Given two points $x, y\in \R^n$ we write
\[
[x, y]:=\{tx+(1-t)y\ : \ 0\leq t\leq 1\},
\]
and we call it the {\it closed segment} between $x$ and $y$ (analogously for $[x, y), (x, y], (x, y)$). For $r>0$ we write $\overline{B}(x, r), B(x, r), C(x, r)$ the closed ball, the open ball and the sphere centered at $x$ with radius $r$, respectively. For $v\neq 0$ in $\R^n$ we write
\[
[x, \infty)_{v}:=\{x+tv\ : \ t\geq 0\}
\]
for the {\it infinite ray starting at $x$ in the direction of $v$}. We write $l_{a, v}:=[a, \infty)_{-v}\cup [a, \infty)_v$ for the entire straight line passing through $a$ in the direction of $v$.
\section{Conics as midsets}\label{conics_as_midsets}
In this section we review the definition of the classical conics as the equidistant set to two circles (possibly degenerating into points or straight lines).
\begin{figure}[H]
  \centering
  \subfigure[Hyperbola]{\includegraphics[width=.3\textwidth]{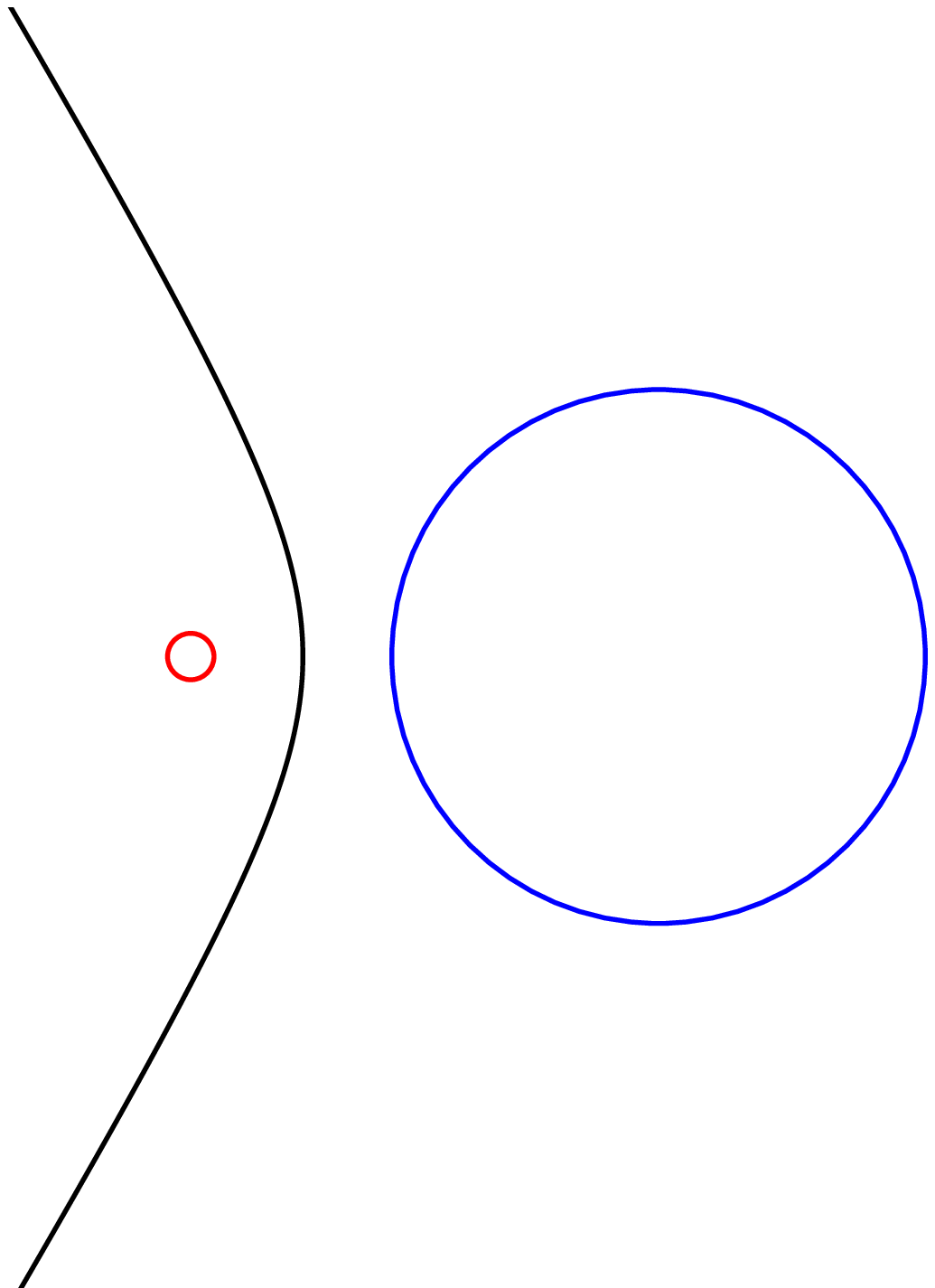}\label{classicHyperbola}}
  \subfigure[Ellipse]{\includegraphics[width=.3\textwidth]{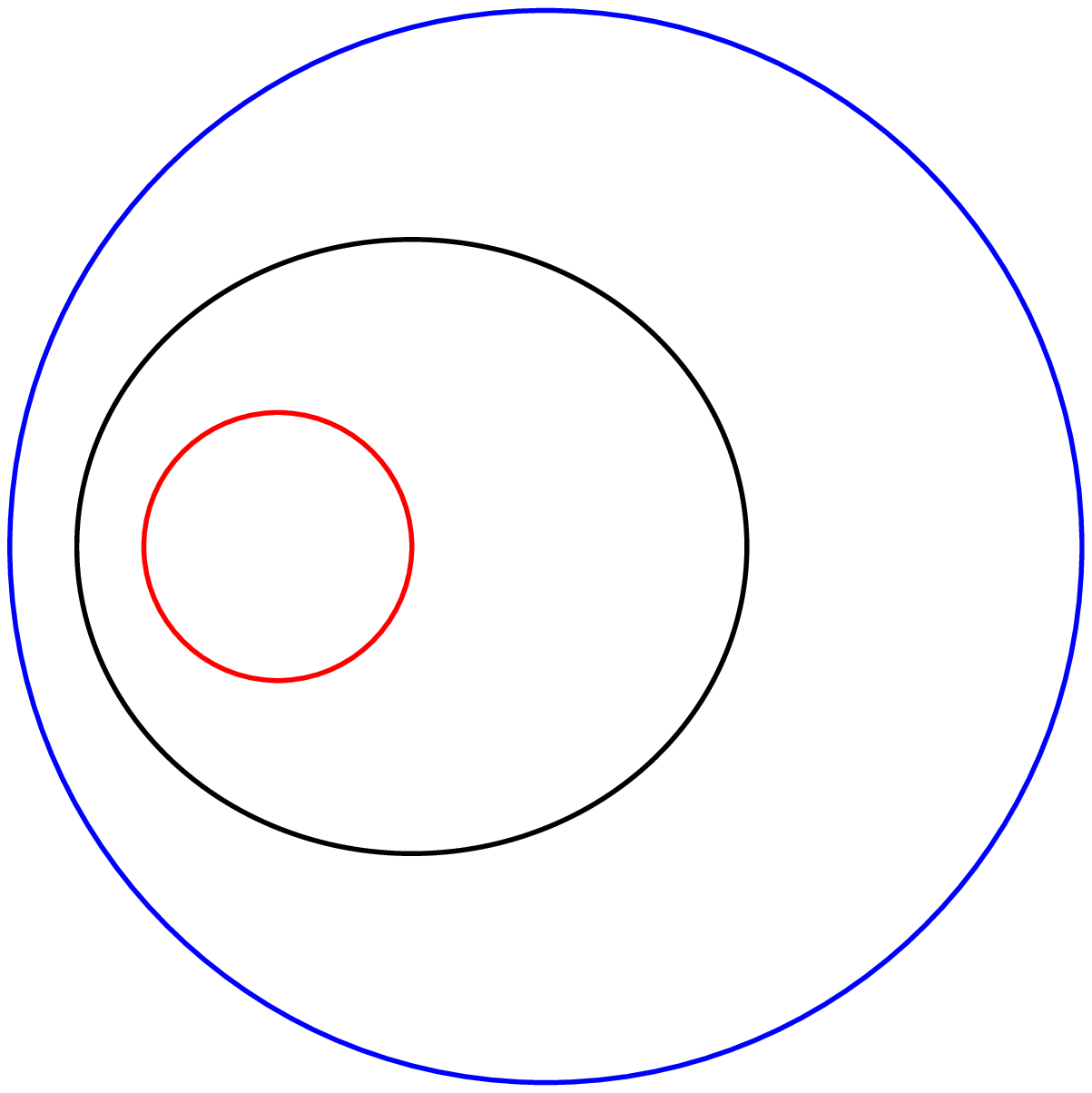}\label{classicEllipse}}
  \subfigure[Parabola]{\includegraphics[width=.3\textwidth]{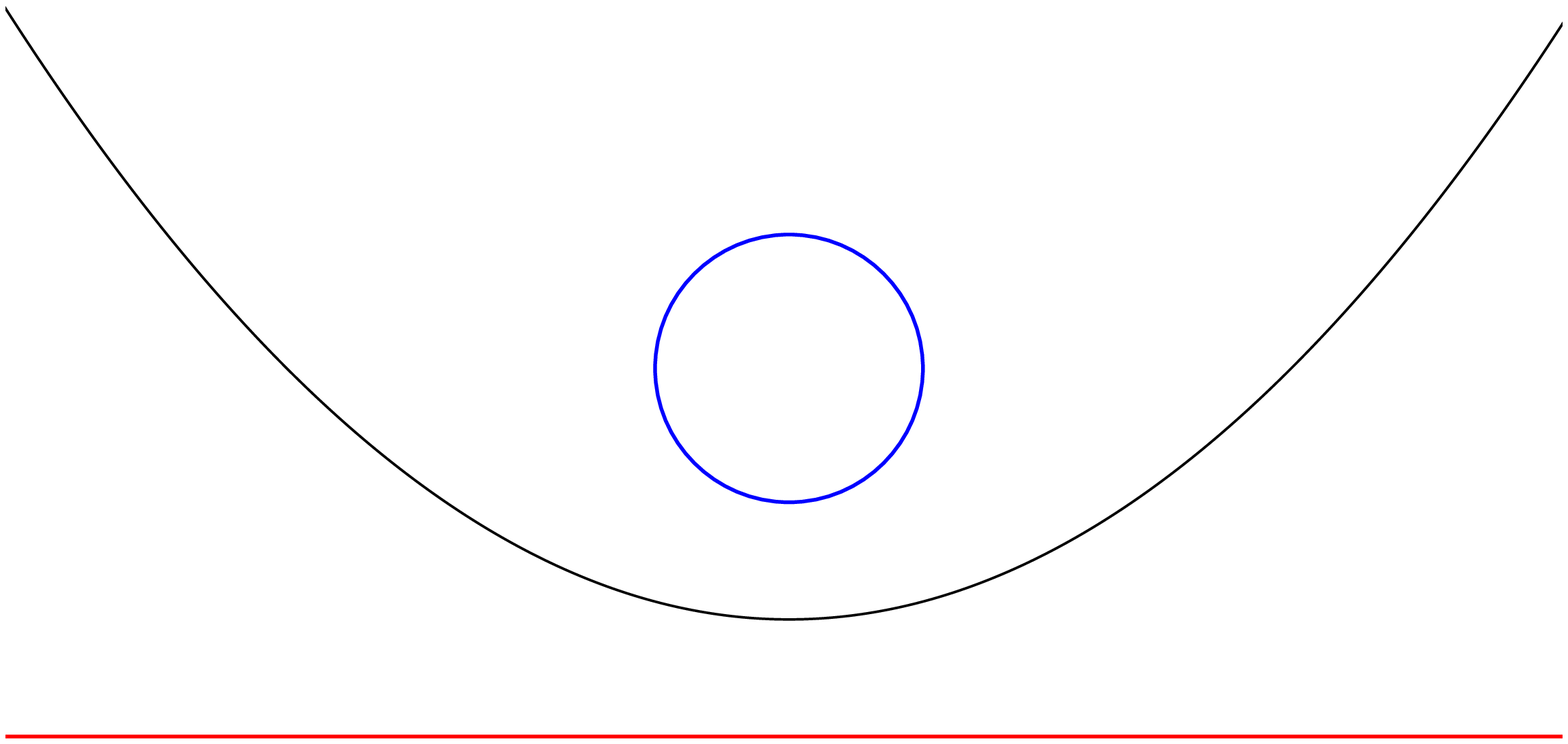} \label{classicParabola}}
  \caption{Classical conics}
  \label{classicConics}
\end{figure}
\paragraph{Hyperbola. } Let $A=C(0, R)$ and $B=C(1, r)$ with $0\leq r, R$ and $R<1-r$ (this implies $A\cap B=\emptyset$).  Using complex notation for points in the plane, the midset $\{A=B\}$ is composed by points $z\in \C$ so that
\begin{eqnarray*}
dist(z, A)&=&dist(z, B)\\
|z|-R&=&|z-1|-r\\
|z-0|-|z-1|&=&R-r.
\end{eqnarray*}
Thus, the midset $\{A=B\}$ is exactly the locus of points $z$ in the plane so that the difference of the distance from $z$ to $0$ and $1$ is constant, that is, the branch of a hyperbola.  In the case $R=r$ we obtain a straight line.
\paragraph{Ellipse. } Consider this time two circles $A=C(0, R)$ and $B=C(1, r)$ with $R>1+r$ (this implies $B$ lies inside $A$). The midset $\{A=B\}$ is composed by points $z\in \C$ so that
\begin{eqnarray*}
dist(z, A)&=&dist(z, B)\\
R-|z|&=&|z-1|-r\\
|z-0|+|z-1|&=&R-r.
\end{eqnarray*}
Thus, the midset $\{A=B\}$ is a ellipse with focal sets $\{0\}$ and $\{1\}$.
\paragraph{Parabola. } The intermediate construction when one of the circles degenerates into a straight line and the other into a point is one of the most classical examples of an equidistant set. Namely, a parabola is the locus of points from where the distances to a fixed point (focus) and to a fixed line (directrix) are equal.
\\

Conversely, we left as an exercise to the reader to show that every ellipse or branch of hyperbola can be constructed as the midset of two conveniently chosen  circles.
\section{Topological properties of midsets}\label{survey}
The first part of this section corresponds to a  survey of some topological properties of midsets in $\R^n$. We mainly concentrate in the articles \cite{WILK75} and \cite{LOVE76}. At sections \ref{continuidad}, \ref{seccion_shadow} we present new results.\\

Since the closure $\overline A$ of a non-empty set $A\subset \R^n$ verifies $dist(x, A)=dist(x, \overline A)$ for every $x\in \R^n$, we can easily conclude that $\{A=B\}=\{\overline A=\overline B\}$, for every $A, B\subset \R^n$.  Hence we can restrict us to consider closed sets as focal sets of midsets. The function
\begin{eqnarray*}
d_A:\R^n&\longrightarrow&\R\\
x&\longmapsto&d_A(x):=dist(x, A)
\end{eqnarray*}
is continuous and
\[
\{A=B\}=d_{A, B}^{-1}(0),
\]
where $d_{A,B}(x):=d_A(x)-d_B(x)$.  We conclude that midsets are always closed sets. Furthermore, a midset is never empty. Indeed,  we can compute the function $d_{A,B}$ over a continuous path joining $A$ and $B$ in order to obtain a zero for $d_{A,B }$ (in fact, it can be shown that every midset is non-empty if and only if   the ambient space is connected). The main theorem in \cite{WILK75} is the following
\begin{theo}[see Theorem 4 in \cite{WILK75}]
If $A$ and $B$ are non-empty connected sets, then $\{A=B\}$ is connected. $\quad_{\blacksquare}$
\end{theo}

 The following is a simple property, that, at least from the point of view of applications to sea frontiers, provides the politically correct fact that there is no region inside the UN's definition of the sea boundary between two disjoint countries.
 \begin{prop}[see \cite{WILK75}]
 Let $A, B\subset \R^n$ be two disjoint non-empty sets. Then the midset $\{A=B\}$ has empty interior.
 \end{prop}
 \noindent{\it Proof. } Let $x\in \{A=B\}$ and let $p_A\neq  p_B$ be foot points in $\mathcal{P}_x(A), \mathcal{P}_x(B)$  respectively.  We claim that for any $\tilde x\in [p_A, x)$ we have
 \begin{equation}\label{clip}
dist(\tilde x, A)<dist(\tilde x, B).
\end{equation}
Indeed, the closed ball $\overline B( x, dist( x, A))$ strictly contains the closed ball $\overline B(\tilde x, dist(\tilde x, A))$. But $dist(\tilde x, p_A)=dist(\tilde x, A)$, which implies that there is no point of $B$ in $\overline{B} (\tilde x, dist(\tilde x, A))$ and then we have (\ref{clip}). The inequality  (\ref{clip}) tell us in particular that $\tilde x\notin \{A=B\}$, and the proposition follows by picking $\tilde x$ as close to $x$ as necessary. $\quad_{\square}$\\
\noindent
\begin{obs}
Notice above that $dist(\tilde x, A)=dist (x, A)-dist(x, \tilde x)$ and (\ref{clip}) can be improved to
\begin{equation}\label{lema_circulo}
dist(\tilde x, B)>dist (x, A)-dist(x, \tilde x).
\end{equation}
\end{obs}
Continuing  with the topological properties of midsets, we concentrate in the case when the focal sets $A, B$ are disjoint compact connected non-empty sets. In that case one has
\begin{theo}\label{elipse_loveland}
\noindent
\begin{itemize}
\item[i) (see \cite{B})] If $A, B\subset \R^2$ then the midset $\{A=B\}$ is a topological $1-$manifold.
\item[ii) (see \cite{LOVE76})] For $A, B\subset \R^n$ with $n\neq 2$ the above result is no longer true in general. However, for every $n$, if $A$ is convex then $\{A=B\}$ is topologically equivalent to an open set of the sphere $\mathcal{S}^{n-1}$. Furthermore, the midset $\{A=B\}$ is homeomorphic to the sphere $\mathcal{S}^{n-1}$ if and only if $A$ is convex and lies in the interior of the convex hull of $B. \quad_{\blacksquare}$
\end{itemize}
\end{theo}
A geometrical objet that is closely related to midsets is the so called $\eps-$boundary of a set $A\in \R^n$
\[
\partial_\eps(A):=\left\{x\in \R^n\ : \ dist(x, A)=\eps \right\}.
\]
In fact one has the relation
\[
\{A=B\}=\bigcup_{\eps\geq 0}\left(\partial_\eps(A)\cap \partial_\eps(B)\right).
\]
These sets have been studied widely in \cite{BROW72}, \cite{FERRY76} and recently in \cite{PIKU03}. A deeper relation between midsets and $\eps-$boundaries is indebted to Loveland:
\begin{theo}[see \cite{LOVE76}]
If $A, B$ are disjoint closes sets of $\R^n$, $A$ is convex and $\eps>0$, then $\{A=B\}$ is homeomorphic to an open subset of $\partial_{\eps}(A). \quad_{\blacksquare}$
\end{theo}
\subsection{Continuity of midsets}\label{continuidad}
 Let $(X, dist_X)$ be a compact metric space. For $A\subset X$ and $\eps>0$ we denote by
 \[
 B(A, \eps):=\{x\in X \ : \ dist_X(x, A)<\eps\}.
 \]
 the $\eps$-neighborhood of $A$. The {\it Hausdorff distance} between two compact sets $K_1, K_2\subset X$ is
 \[
 dist_{\mathcal{H}}(K_1, K_2):=\inf \{\eps>0 \ : \ K_1\subset B(K_2, \eps)\ \textrm{and}\ K_2\subset B(K_1, \eps)\}.
 \]
This distance defines a topology on the space $\mathcal{K}(X)$, of compact subsets of $X$. With this topology  $\mathcal{K}(X)$ is itself a compact space (see for instance \cite{MONTZIP74}). Given a convergent sequence $A_n\in \mathcal{K}(X)$, the  Hausdorff limit is characterized  as the set of points that are limits of sequences $x_n\in A_n$.
\\

In general the equidistant sets are closed but not necessarily bounded sets.  In order to treat with compact sets and use the Hausdorff topology, we are going to consider restrictions of equidistant sets to a large enough ball containing both focal sets.  Let $R$ be a large positive number and $A, B$ be compact sets so that $A\cup B\subset B(0, R)$. We write
\[
\{A=B\}_R:=\{A=B\}\cap \overline B(0, R).
\]
We are interested in the continuity of the application
\begin{eqnarray*}
\mathcal{M}id_R:\mathcal{K}(\overline B(0,R))\times \mathcal{K}(\overline B(0,R))&\longrightarrow& \mathcal{K}(\overline B(0,R))\\
(A, B)&\longmapsto&\{A=B\}_R.
\end{eqnarray*}
\begin{theo}\label{hausdorff}
  If $A\cap B=\emptyset$ then $(A, B)$ is a continuity point of $\mathcal{M}id_R$.
\end{theo}
\noindent{\it Proof. } Let $\{A_n\}_{n\in \N}, \{B_n\}_{n\in \N}$ be two sequences in $\mathcal{K}(\overline B(0,R))$ so that
\[
A_n\to A\quad \textrm{and}\quad B_n\to B.
\]
Define $E_n:=\{A_n=B_n\}_R\in \mathcal{K}(\overline B(0,R))$. A compactness argument  allows to assume that there exists $E\in \mathcal{K}(\overline B(0,R))$ so that $E_n\to E$. We affirm  that $\{A=B\}_R=E$.
\begin{itemize}
\item Let $e\in E$. There exists sequences $e_n\in E_n$, $a_n \in A_n$, $b_n\in B_n$ and two points $a\in A$, $b\in B$ so that
\begin{equation}\label{tres_igual}
dist(e_n, a_n)=dist(e_n, A_n)=dist(e_n, B_n)=dist(e_n, b_n)
\end{equation}
with $e_n\to e$, $a_n\to a$ and $b_n\to b$.  We claim that $dist(e, A)=dist(e, a)$. Assume otherwise that there is a point $\tilde a\in A$ so that  $
dist(e, \tilde a)<dist(e, a)$.
There exists a sequence $\tilde a_n\in A_n$ with $\tilde a_n\to \tilde a$. But (\ref{tres_igual}) implies
\[
dist(e_n, a_n)\leq dist(e_n, \tilde a_n),
\]
that leads $dist(e, a)\leq dist(e, \tilde a)$. In a similar way one shows $dist(e, B)=dist(e, b)$. Taking limit in (\ref{tres_igual}) we get $dist(e, A)=dist(e, B)$ and then $E\subset \{A=B\}_R$.
\item Let $m\in \{A=B\}_R$ and $a_n\in A_n$, $b_n\in B_n$ verifying
\begin{eqnarray*}
dist(m, A_n)&=&dist(m, a_n),\\
 dist(m, B_n)&=&dist(m, b_n).
\end{eqnarray*}
Passing to a subsequence if necessary there exist $a\in A$, $b\in B$ so that $a_n\to a$ and $b_n\to b$. Then
\begin{eqnarray}
dist(m, a_n)=dist(m, A_n)&\to & dist(m, a), \\
dist(m, b_n)=dist(m, B_n)&\to & dist(m, b).
\end{eqnarray}
From this one has
\[
\lim_{n\to \infty}dist(m, A_n)-dist(m, B_n)=0.
\]
Passing to a subsequence (or interchanging the roles of $A_n$ and $B_n$) we can assume that $dist(m, A_n)-dist(m, B_n)$ increases to zero. Let $t\geq 0$. We define $m_t\in [m, b]$ so that $dist (m, m_t)=t$.  Define $f_n(t)=dist(m_t, A_n)-dist(m_t, B_n)$.  Let $\eps>0$. We claim that there exists $\tilde n\in \N$ so that $f_{n}(\eps)>0$ for every $n\geq \tilde n$. Indeed, we know that
\begin{eqnarray*}
dist(m_\eps, B)&=&dist(m, B)-\eps ,  \\
dist(m_\eps, A)&>& dist(m, A)-\eps.
\end{eqnarray*}
Notice that the second inequality above follows from (\ref{lema_circulo}) (here we use the hypothesis $A\cap B=\emptyset$). From these we obtain $dist(m_\eps, A)-dist(m_\eps, B)>0$. Since $f_n(\eps)\to dist(m_\eps, A)-dist(m_eps, B)>0$ our claim holds.  Using that $f_n(0)\leq 0$ for every $n\geq \tilde n$ we can pick $m_n\in [m, m_\eps]$ so that $f_n(m_n)=0$, that is $m_n\in \{A_n=B_n\}$. This construction holds for every $\eps>0$, and then a diagonal sequence argument allows to construct a sequence $m_n\in E_n$ with $m_n\to m$. That is, $m\in \lim E_n=E$, and finally $\{A=B\}_R\subset E. \quad_{\square}$
\end{itemize}
\section{Midsets as generalized conics}
In section \ref{conics_as_midsets} we have seen how the classical conics can be realized as equidistant sets with circular focal sets. In this section we want to interpret equidistant sets as na\-tu\-ral generalizations of conics when admitting  focal sets that are more complicated than circles. We concentrate in recovering geometric properties from conics to more general midsets. \\
\begin{figure}[H]
  \centering
  \subfigure[Hyperbola]{\includegraphics[width=.3\textwidth]{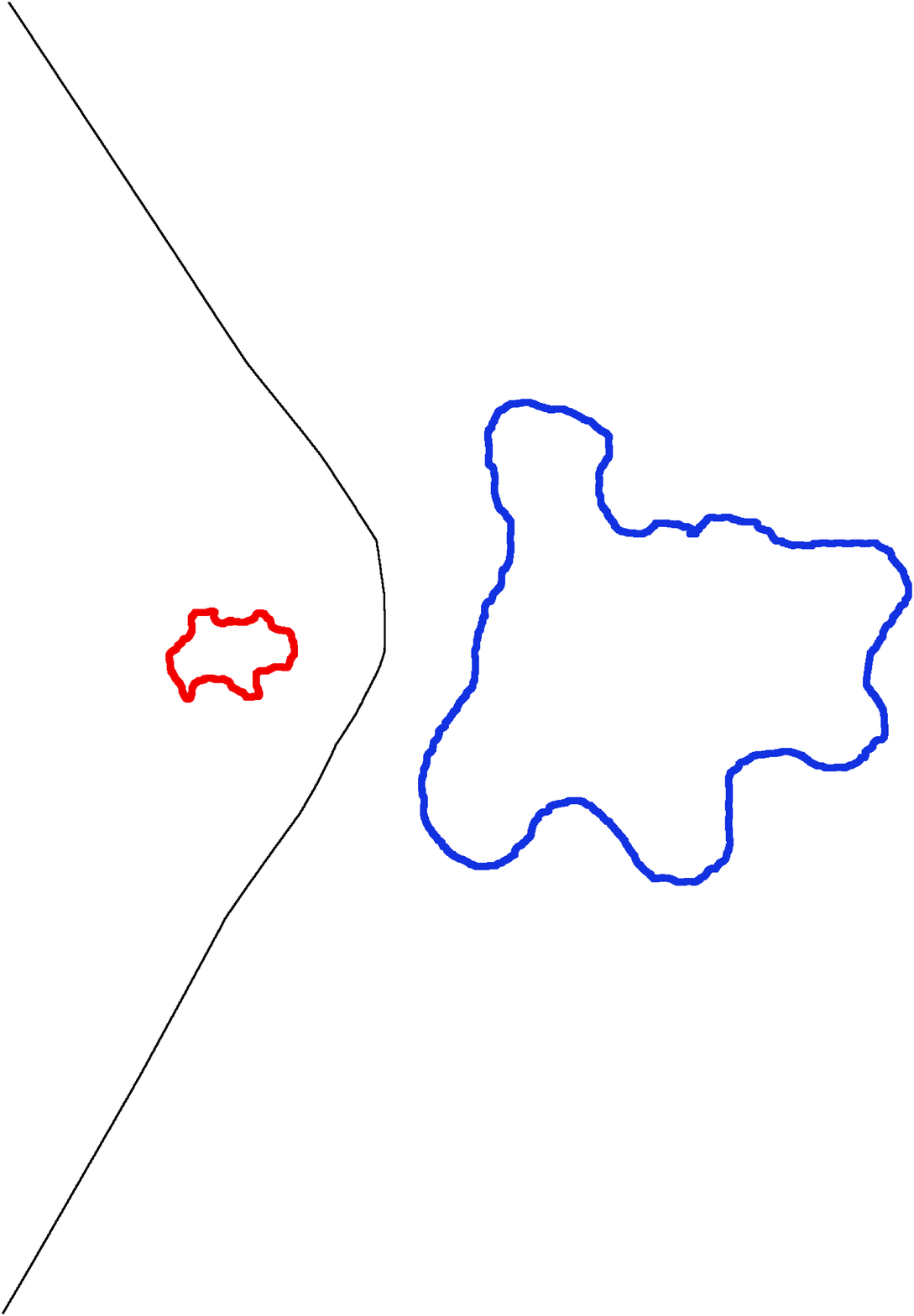} \label{generalizedHyperbola}}
    \subfigure[Ellipse]{\includegraphics[width=.3\textwidth]{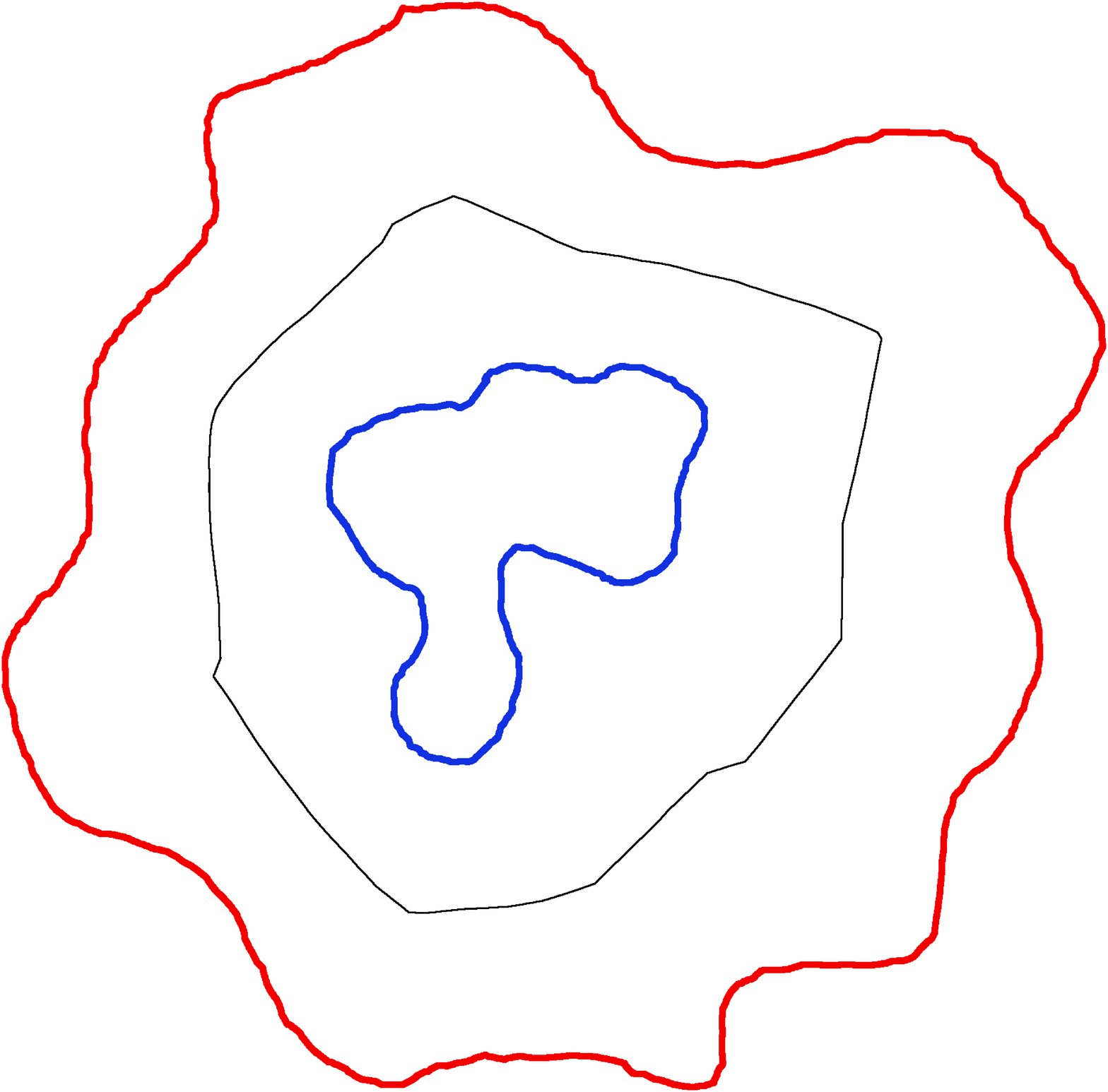}\label{generalizedEllipse}}
  \subfigure[Parabola]{\includegraphics[width=.3\textwidth]{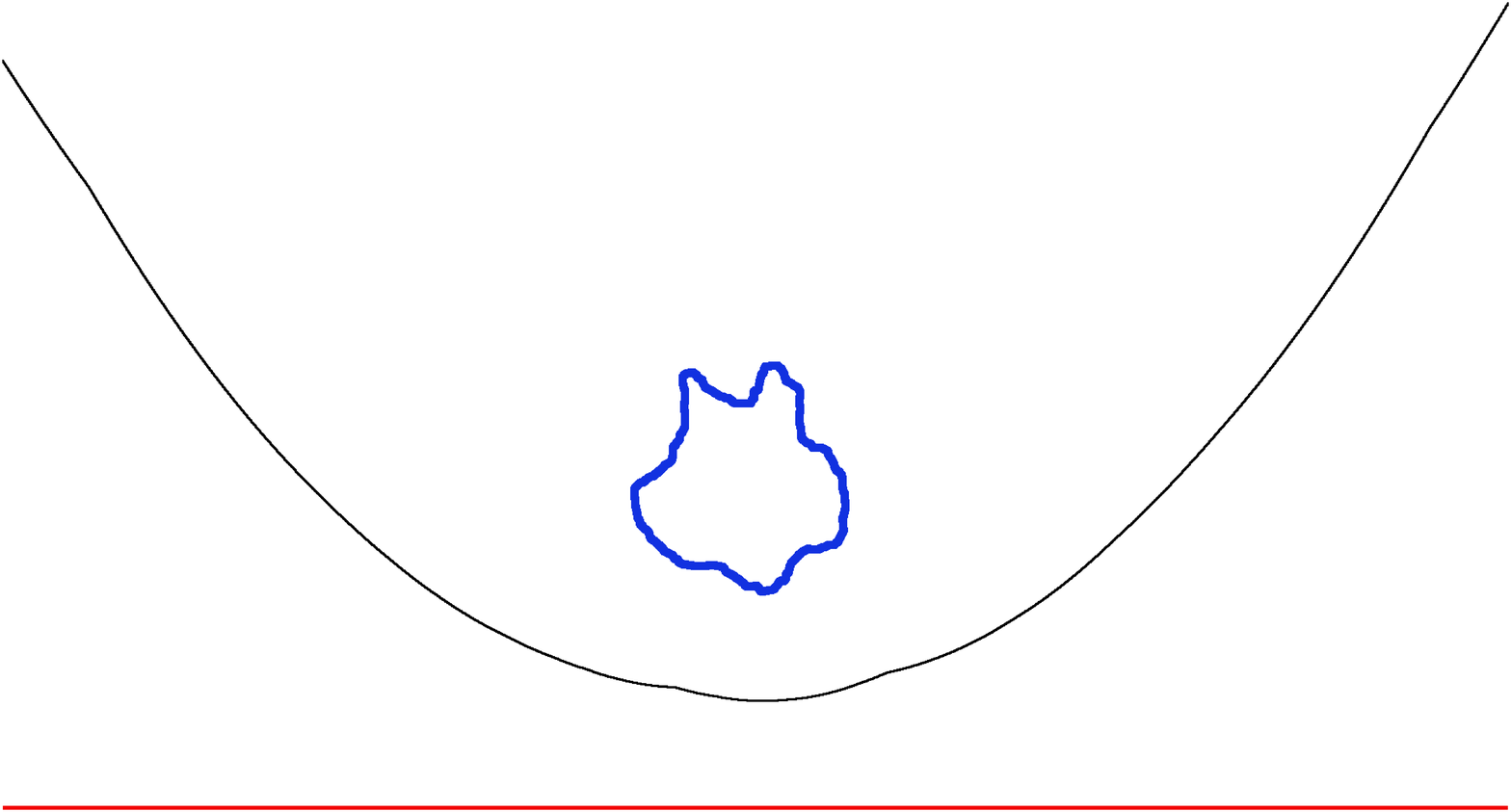}\label{generalizedParabola}}
  \caption{Generalized conics}
  \label{generalizedConics}
\end{figure}
\paragraph{Generalized hyperbolas. } In section \ref{conics_as_midsets} we have seen that a branch of a hyperbola can be realized as the midset of two disjoint circles. In this section we show that replacing these two discs by two {\it disjoint enough} compact connected sets we recover a midset that asymptotically resembles a branch of a hyperbola. Indeed, we show that far enough from these focal sets  the midset consists in two disjoint continuous curves that go to infinity  approaching asymptotically to two different directions in the plane. This is the content of Theorem \ref{dos_colas} bellow.  \\

We require some additional definitions and notation. Let $\vec r=[a, \infty)_{v}$ be a ray starting at $a\in \R^2$ in direction $v\in \R^2$, with $\|v\|=1$. Pick $v^{\perp}$ so that $\{v, v^{\perp}\}$ is a positive orthonormal basis for $\R^2$. For $\eps>0$ we define the {\it tube} of width $\eps$ around $\vec r$ as
\[
tub_{\eps}(\vec r):=\left\{a+tv+sv^{\perp}\  : \ t\geq 0, |s|\leq \eps\right\}
\]
We say that a set $M$ has  an {\it asymptotic end in the direction of $\vec r$} if there exists $\eps>0$ so that the set $M_{\eps, \vec r}=M\cap tub_{\eps}(\vec r)$ verifies the two following conditions:
\begin{itemize}
\item[$i)$ ] The orthogonal projection from $M_{\eps, \vec r}$ to $\vec r$ is a bijection.
\item[$ii)$ ] If we write $M_{\eps, \vec r}$ using the parameters $(t, s)$ of the tube $tub_{\eps}(\vec r)$, then the point $(i)$ above yields a function
\begin{eqnarray*}
s:[0, \infty)&\longrightarrow&[-\eps, \eps]\\
t&\longmapsto&s(t)
\end{eqnarray*}
in such a way that $M_{\eps, \vec r}$ coincides with the graph of $s$. The second requirement  is that
\[
\lim_{t\to \infty}s(t)=0.
\]
\end{itemize}

\begin{obs}
\noindent
\begin{itemize}\item[$i) $] Notice that the function $s$ defined above is continuous since its graph is a closed set.
\item[$ii)$ ] The reader can notice that every ray $\vec p\subset \vec r$ induces an asymptotic end just by considering the suitable restriction. Even though one can formalize properly using an equivalence relation, we are going to consider all these ends as the same one.
\end{itemize}
\end{obs}
Let $K\subset \R^2$ be a compact set. We say that the straight line $l=l_{b, w}$ is a {\it supporting line} for $K$ if $l\cap K\neq \emptyset$ and $K$ is located entirely in one of the two half-planes defined by $l$.
We say that $x\in l\cap K$ is a {\it right extreme point} if $l\cap K$ is contained in $[b, \infty)_{-w}$ (analogously we define an {\it left extreme point}). A supporting line has  always both types of extreme points, and they coincide if and only if the intersection $l\cap K$ contains only one point.
\begin{lema}\label{lrma9}
Let $\eps>0$. Assume that $K\subset \{(x, y) \ | \ x \leq \eps\ , \ y \leq 0 \}$. For $h>0$ we define $f_h(x)=dist\left((x, h), K\right)$. The function  $f_h$ is strictly  increasing for $x\geq \eps$.
\end{lema}
\noindent{\it Proof. } Let $x_2>x_1 \geq \eps$ and let $p_2\in \mathcal{P}_{(x_2, h)}(K)$ be a foot point.   We have
\[
f_h(x_1)\leq dist ((x_1, h), p_2)<dist((x_2, h), p_2)=f_h(x_2).
\]
Indeed, the first inequality comes from the definition of $f_h(x_1)$ and the strict inequality is due to $x_2>x_1. \quad_{\square}$ \\

In what follows we are going to consider two compact sets $A, B$ and a common supporting line $l$ so that both sets are located in the same half-plane determined by $l$. For simplicity we assume that $l$ is the real line and $(-1,0)$ is the right extreme point of $A$ and $(1, 0)$ is the left extreme point of $B$.  Let $\eps>0$ small enough. We assume that
\begin{eqnarray}\label{drop_A}
A&\subset& \{(x, y)\ : \ x\leq -1+\eps \ , \ y\leq 0\}, \\ \label{drop_B}
B&\subset& \{(x, y)\ : \ x\geq 1-\eps \ , \ y\leq 0\}.
\end{eqnarray}
\begin{lema}\label{estrella_sol}
Under the above hypotheses, for every $h>0$  there exists a unique $x(h)\in [-1, 1]$ so that
\begin{equation}\label{ecu_de_lema10}
dist((x(h), h), A)=dist((x(h), h), B).
\end{equation}
Moreover, $x(h)$ belongs to $(-\eps, \eps)$.
\end{lema}
\noindent{\it Proof. } Since $(-1,0)\in A$, we know that for every $(x, y)\in \{x\leq -\eps\ ,\ y\geq 0 \}$ one has
\[
dist((x, y), A)<dist((x, y), B).
\]
Similarly we obtain that for every $(x, y)\in \{x\geq \eps\ ,\ y\geq 0 \}$ one has
\[
dist((x, y), A)>dist((x, y), B).
\]
Then the continuity of the function $f_h$ defined in Lemma \ref{lrma9} gives at least one point $x(h)\in (\-eps, \eps)$ satisfying the equality (\ref{ecu_de_lema10}). Applying the conclusion of Lemma \ref{lrma9} we see that the function
\[
x\longmapsto dist((x, h), A)-dist((x, h), B)
\]
 is strictly increasing for $x\in [-1+\eps, 1-\eps]$. We then deduce the unicity of $x(h)$ as required. $\quad_{\square}$
\\

We apply the above Lemma in order to characterize asymptotically the midset of two focal sets with a common supporting line. Notice that in the hypotheses of the next Proposition we drop conditions (\ref{drop_A}) and (\ref{drop_B}).
\begin{prop}\label{prop_ONCE}
Consider two disjoint compact sets $A, B$ and a common supporting line $l$ so that both sets are located in the same half-plane determined by $l$. For simplicity we assume that $l$ is the real line and $(-1,0)$ is the right extreme point of $A$ and $(1, 0)$ is the left extreme point of $B$. For every $\eps>0$ there exists $\tilde h=\tilde h(\eps)>0$ such that for every $h>\tilde h$ there exists $x(h)\in (-\eps, \eps)$ so that the following holds
\[
\{A=B\}\cap \{(x, h)\ : \ x\in [-1,1]\}=\{(x(h),h)\}.
\]
\end{prop}
\noindent{\it Proof. } In order  to apply Lemma \ref{estrella_sol}, we need to show that we can recover conditions  (\ref{drop_A}), (\ref{drop_B}). Since $(-1,0)$ is in $A$, for every $h>0$ the foot points $\mathcal{P}_{(0, h)}(A)$ belongs to  the closed ball $D_h$ centered at $(0, h)$ and passing through $(-1,0)$  (the same happens for $\mathcal{P}_{(0, h)}(B)$, for the same ball $D_h$ since also passes trough $(1,0)$). In other words one has \[
\mathcal{P}_{(0, h)}(A)\cup \mathcal{P}_{(0, h)}(B)\subset D_h\cap \{(x, y)\ | \ y\leq 0\}.
\]
We define
\[
A_{h}:=D_h\cap A\quad, \quad B_{h}:=D_h\cap B.
\]
With these definitions it is clear that
\begin{eqnarray*}
dist((0, h), A)&=&dist((0, h), A_h), \\
dist((0, h), B)&=&dist((0, h), B_h).
\end{eqnarray*}
We claim that for every $\eps >0$ there exists $\tilde h>0$ so that for every $h>\tilde h$ one has
\begin{eqnarray*}
A_h&\subset& \{(x, y)\ | \ x\leq -1+\eps \ , \ y\leq 0\}, \\
B_h&\subset& \{(x, y)\ | \ x\geq 1-\eps \ , \ y\leq 0\}.
\end{eqnarray*}
Assume in the contrary that there exists $\tilde \eps>0$ and a sequence $(x_n, y_n)\in A_{n}$ with  $x_n>-1+\tilde \eps$. Notice that $(x_n, y_n)\in D_{n}$, then we have
\[
-y_n+n\leq \sqrt{n^2+1}.
\]
This and the classical undergraduate limit $\lim_{n\to \infty} \sqrt{n^2+1}-n=0$  implies that $y_n\to 0$. Since $A$ is a compact set, there exists a subsequence $(x_n, y_n)$ converging to a point $(\tilde x, 0)\in A$, with $\tilde x\geq -1+\tilde \eps>-1$. This contradicts the fact that $(-1,0) $ is the right extreme point of $A$. We then apply Lemma \ref{estrella_sol} in order to find $x(h)\in (-\eps, \eps)$ in the midset $\{A=B\}$. It is easy to see that for fixed $\eps>0$ and $h$ large enough we have $(-\eps, h)$ is closer to $A$ and $(\eps, h)$ is closer to $B$, concluding thus the proof. $\quad_{\square}$\\

Given two disjoint non-empty compact connected sets $A, B$, we want to discuss about the existence of a common supporting line leaving both sets in the same half-plane. For this we need to remember the concept of convex hull $ch(K)$ of a compact set $K\subset \R^2$ defined as the smallest convex set containing $K$. The convex hull $ch(K)$ is a convex compact set. Given two disjoint compact convex sets $\mathcal{A}, \mathcal{B}\subset \R^2$, it is an interesting exercise to show that there exist four common supporting lines.  Two of them are called {\it interior common tangents} and each one  leaves the sets $\mathcal{A}, \mathcal{B}$ into a different half-plane. The remaining two supporting lines are called {\it exterior common tangents} and each one leaves both sets into the same half-plane. \\
\begin{figure}[H]\nonumber
\centering
 \includegraphics[width=11cm]{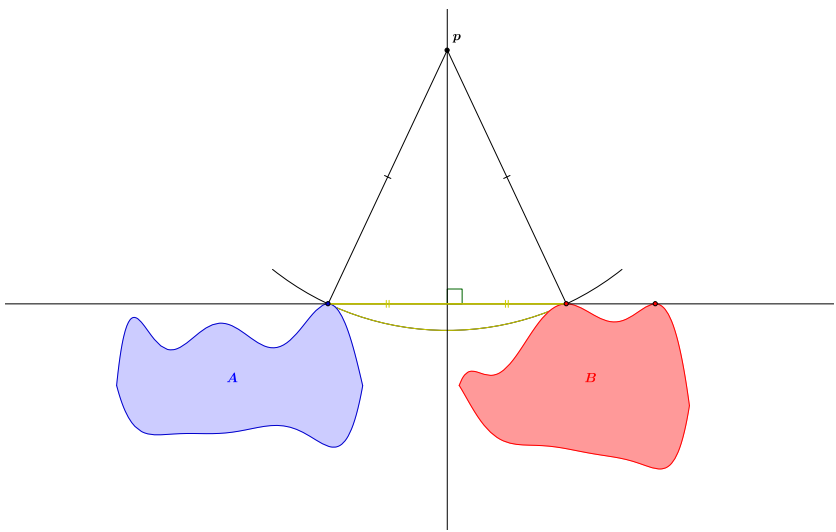}
  \caption{Foot points to $p$ lie inside the small chordal region. }
\end{figure}
 Two disjoint non-empty compact connected sets $A, B$ are called {\it $ch$-disjoint } if $ch(A)\cap ch(B)=\emptyset$. It is easy to see that supporting lines and common supporting lines of $ch(A), ch(B)$ are also supporting lines and common supporting lines of $A, B$ respectively. The above discussion directly yields
 \begin{lema}\label{lema_DOCE}
 Two non-empty compact connected sets $A, B$ that are {\it $ch$-disjoint } have two distinct common supporting lines each of one leaves both sets $A, B$ into the same half-plane . $\quad_{\square}$
  \end{lema}
  Now we can state the main theorem of this section
  \begin{theo}[Generalized hyperbola]\label{dos_colas}
  Let $A, B$ be two non-empty compact connected sets that are {\it $ch$-disjoint}. There exists $R>0$ and two disjoint rays $\vec r_1, \vec r_2$ so that
  \[
  \{A=B\}\cap B(0, R)^c
  \]
  consists exactly of two asymptotic ends in the direction of $\vec r_1$ and $\vec r_2$ respectively.
  \end{theo}
  \noindent{\it Proof. } The existence of the two different asymptotic ends is due to Proposition \ref{prop_ONCE} and Lemma \ref{lema_DOCE}. The remaining part of the proof consists in to show that there is no other piece of the midset going to infinity. This can be directly deduced from Theorem \ref{elipse_loveland} part $(i)$ which ensures that the midset $\{A=B\}$ is homeomorphic to the real line. We also present a self-contained proof: assume that $\vec r_1, \vec r_2$ are not parallel and suppose that there exists a sequence $x_n\in \{A=B\}$ with $|x_n|>n$. Let $x_n=|x_n|e^{i\theta_n}$ be the complex notation for $x_n$. Taking a subsequence if needed, we can  assume that there exists $\tilde \theta\in [0, 2\pi]$ so that $\theta_n\to \tilde \theta$. Let $l_{\tilde \theta^{\perp}}$ be a supporting line for $A\cup B$ that is orthogonal to the direction $\tilde \theta$ and so that $A\cup B$ and $\{x_n\}_{n\in \N}$ are located in different half-planes. We claim that $l_{\tilde \theta^{\perp}}$ is a common supporting line for $A$ and $B$. Indeed, assume for instance  that $l_{\tilde \theta^{\perp}}\cap A=\emptyset$. In this case it is easy to see that for $n$ large enough we should have $dist(x_n, A)>dist(x_n, B)$, that is impossible since $x_n$ belongs to $\{A=B\}$.  Hence, $\tilde \theta$ coincides with the direction of $\vec r_1$ (or $\vec r_2$) and necessarily one deduces that $\{x_n\}$ is a subset of the union of the two asymptotic ends. The case when $\vec r_1, \vec r_2$ are parallel can be easily treated by considering a slight perturbation  $A_\eps$ of $A$ in the Hausdorff topology in such a way that the above lines can be applied to $A_\eps, B$. We conclude using Theorem \ref{hausdorff}. $\quad_{\square}$
  \begin{obs}
  For simplicity we stated this theorem for $ch$-disjoint sets, even though it holds for every pair of compact connected disjoint sets having two supporting lines each of one leaves both sets in the same half-plane.
  \end{obs}
\paragraph{Generalized ellipses. } In this case we have not much to say. However, the part $(ii)$ of Theorem \ref{elipse_loveland} serves to recognize  some topological reminiscences of ellipses when the focal sets of the midset are a convex compact set {\it inside} a  compact set.
 \paragraph{Generalized parabolas. } A remarkable geometric property of parabolas is that they are {\it strictly convex} in the sense that for any supporting line the asymptotic behavior  at infinity consists in to become more and more separated from the supporting line. In other words, the parabola can be seen as the graph of a continuous function over the supporting line (a tangent) so that the values of this function tend to infinity with the parameter of the line (check for instance the parabola $y=x^2$ and see how the derivatives grow to infinity). We are not going to give a definition for generalized parabolas . Instead we want to say that midsets sharing some pro\-per\-ties like {\it strict convexity} should be considered as some kind of generalization for parabolas.

 Along the lines of the generalized hyperbolas treated in the previous paragraphs, we want to consider a midset defined by a compact connected focus $A$ (in the place of the classical focus point) and some disjoint unbounded closed set $B$ playing the role of the directrix.  We also need to require some additional properties like: the $ch(B)$ does not intersect $A$ (in order to obtain an unbounded midset); there is no common   supporting line for $A$ and $B$ (in order to avoid the existence of an asymptotic ray), etc. For simplicity we are going to keep $B$ as a straight line, even though the reader will be able to treat with more general situations.
 \begin{prop}
 Let $A\subset \R^2$ be a non-empty connected compact set and $B$ be a disjoint straight line. There exists $\overline R>0$ so that for every $R\geq \overline R$ and every supporting line $l$ for $\{A=B\}\cap B(0, R)^c$ one has
 \[
 \lim_{s\to \infty}dist \left(l, \{A=B\}\cap B(0, s)^c\right)=\infty.
 \]
 \end{prop}
\noindent{\it Sketch of the proof. } For this special case where $B$ is a straight line the proof can be easily obtained from the fact that the midset $\{A=B\}$ actually is the graph of a continuous function over $B$. We left as an exercise to the reader to show that this function grows faster than any linear map.

We want to outline a proof that fits to more general situations. Assume for simplicity that $B$ is the real line. The idea  is to truncate $B$ and consider the midset $\{A=B_R\}$, where $B_R=[-R, R]\subset \R$. As seen before, this is a generalized hyperbola that is asymptotic (lets say, to the right) to a ray $\vec r_R$ that is perpendicular trough the midpoint of a segment $[ a_R, (R, 0)]$, for some point $a_R\in A$. Since $A$ is compact, the slope of $\vec r_R$ grows to infinity with $R$ and the reader can easily complete the details. $\quad_{\square}$
\paragraph{Concluding remarks. } The final part  of this  paper  concentrates in to show that equidistant sets looks simple at least from the asymptotic point of view. As a future line of research we suggest to explore how complicated can actually be an equidistant set.  As commented in the Introduction, this question should become crucial since equidistant sets are meant to be used as  region boundaries for many real life situations. For instance,  a sea delimitation needs to have some {\it physical} properties in order to make possible  its role for the real life.  In the measure the equidistant sets admit more a more intricate structures, it become more and more difficult to consider these sets as viable frontiers.  For example, we suggest the following:
\\

\noindent{\it Question: } Does there exist an equidistant set in the plane, with connected disjoint focal sets, having Hausdorff dimension larger than $1$? What about other notions of dimension? How the dimension of the equidistant set depends on the dimension of the focal sets?
\\

\noindent{\it Question: } To characterize all closed sets of $\R^2$ that can be realized as the equidistant set of two connected disjoint closed sets.

\appendix
\section{The Shadowing property}\label{seccion_shadow}
Given two non-empty closed disjoint sets $A, B$ and $\eps>0$, we define the set of $\eps-$equidistant points to $A$ and $B$ as
\[
\{|A-B|<\eps\}:=\left\{x\in \R^n \ : \ |dist(x, A)-dist(x, B)|<\eps\right\}.
\]
This notion is crucial when we deal with computer simulations. Recall that finding an equidistant point is equivalent to find a zero of a continuous function. In the case of computer simulations, this function is no longer continuous since it is evaluated in pixels (a discrete set). In fact, this function in general may have no zero at all. Then, in order to draw a good picture of the equidistant set we need to  check for points (pixels) so that the difference between the distances to the focal sets is small enough to guarantee that inside a small neighborhood there is  a zero for the continuos function   that defines the midset. In conclusion, we look for a set $\{|A-B|<\eps\}$ for some positive $\eps$ that depends on the screen resolution, computer capabilities, etc.
As we will see, the theorem we present here requires a very specific configuration of the focal sets. Nevertheless, the reader should notice that the result can be applied to more general situations.\\

Let  $x\notin A\cup B$. We say that {\it $x$ see $A$ separated from $B$ by an angle $\alpha$} if there exists two supporting lines $l_A, l_B$ passing through $x$ so that
\begin{enumerate}
\item $l_A$ is a supporting line for $A$, and $B$ lies in a different half-plane than $A$.
\item $l_B$ is a supporting line for $B$, and $A$ lies in a different half-plane than $B$.
\item The angle formed at $x$ by $l_A$ and $l_B$ is $\alpha$.
\end{enumerate}
\begin{figure}[H]\label{figura_rara}
\centering
  \includegraphics[width=11cm]{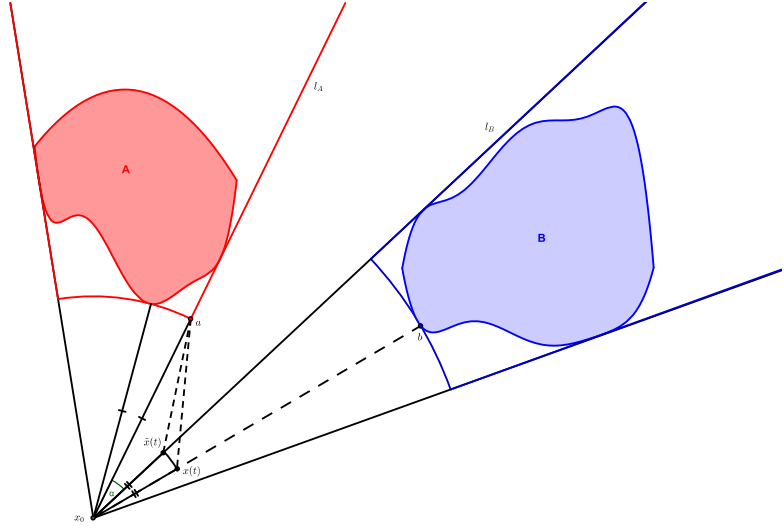}
 \caption{Construction for the proof of Theorem \ref{shadow}}
\end{figure}
\begin{theo}[Shadowing property]\label{shadow}
Let $A, B$ be two disjoint non-empty closed sets . Let $\eps>0$ and $x_0\in \{|A-B|<\eps\}$ so that $x_0$ see $A$ separated from $B$ by an angle $\alpha$. Then there exists $x_1\in \{A=B\}$ verifying
\[
dist(x_1, x_0)<\frac{\eps}{2}\left(\frac{\eps+2d}{\eps+d-d\cos \alpha}\right),
\]
where $d=\min \{dist(x_0, A), dist(x_0, B)\}$.
\end{theo}
\noindent{\it Proof. } Consider $f(x):=d_{B, A}(x)=dist(x, B)-dist(x, A)$ and assume $d=dist(x, A)$, that is
\[
0<f(x)<\eps.
\]
We look for $x_1$ so that $f(x_1)=0$. Let $b\in \mathcal{P}_{x_0}(B)$. We write $x(t)$ as the point on $[x_0, b]$ so that $dist(x_0, x(t))=t$. Finally we write
\[
g(t):=f(x(t))=dist(x(t), B)-dist(x(t), A).
\]
Since $f(x_0)=g(0)$ we have $0<g(0)<\eps$ and $g(\overline d)<0$ where $\overline d=dist(x_0, B)$.  Although we know that there exists $\overline t\in (0, \overline d)$ so that $g(\overline t)=0$, the function $g$ is not differentiable in general and we can not estimate  directly the size of $\overline t$. We are going to construct an upper bound for $g$ in order to get a good estimate. Let $a$ be the intersection of the circle centered at $x_0$ and radius $d$ with $l_A$ as pointed in the figure 4. For every $t$ we have
\begin{equation}\label{primera_estrella}
dist(x(t), a)\leq dist(x(t), A).
\end{equation}
Define $\tilde x(t)\in l_B$ so that $dist(x_0, \tilde x(t))=t$ (see figure 4.). Thus we have
\begin{equation}\label{segunda_estrella}
dist(\tilde x(t), a)\leq dist(x(t), A).
\end{equation}
The left term above can be explicitly computed using elementary euclidean geometry:
\begin{equation}\label{tercera_estrella}
dist(\tilde x(t), a)^2=d^2+t^2-2dt\cos (\alpha).
\end{equation}
Moreover, we know that
\begin{eqnarray}\nonumber
dist(x(t), B)&=&dist(x_0, B)-t, \\ \label{cuarta_estrella}
&<&d+\eps -t.
\end{eqnarray}
Using (\ref{primera_estrella}, \ref{segunda_estrella}, \ref{tercera_estrella}, \ref{cuarta_estrella}) one gets (and define $\hat g$ by)
\begin{equation}\label{quinta_estrella}
g(t)<d+\eps-t -\sqrt{d^2+t^2-2dt\cos (\alpha)}:=\hat g(t).
\end{equation}
Notice that $\hat g(0)=\eps$ and
\[
\hat t=\frac{\eps}{2}\left(\frac{\eps+2d}{\eps+d-d\cos \alpha}\right)
\]
verifies $\hat g(\hat t)=0$. Finally the inequality  (\ref{quinta_estrella}) help us to find a point $\bar t\in (0, \hat t)$ so that $f(\bar t)=0. \quad_{\square}$
\\
\begin{small}

\noindent{\bf Acknowledgments.} Both authors were funded by the  the Fondecyt Grant 11090003.

\end{small}

\begin{footnotesize}

\vspace{0.25cm}

\noindent{Mario Ponce}

\noindent{Facultad de Matem\'aticas, Universidad Cat\'olica de Chile}

\noindent{Casilla 306, Santiago 22, Chile}

\noindent{E-mail: mponcea@mat.puc.cl}
\\

\noindent{Patricio Santib\'a\~nez}

\noindent{Facultad de Matem\'aticas, Universidad Cat\'olica de Chile}

\noindent{Casilla 306, Santiago 22, Chile}

\noindent{E-mail: patriciosantib@gmail.com}

\end{footnotesize}


\end{document}